\documentclass{article}
\usepackage{indentfirst}

\usepackage[english]{babel}

\usepackage[letterpaper,top=2cm,bottom=2cm,left=3cm,right=3cm,marginparwidth=1.75cm]{geometry}

\usepackage{amsfonts,amssymb,amsmath,amsthm, multicol, caption, subcaption, floatrow,tikz}
\usepackage{graphicx}
\usepackage[colorlinks=true, allcolors=blue]{hyperref}

\newtheorem{thm}{Theorem}[section]

\newtheorem{lem}[thm]{Lemma}

\newtheorem{rmk}[thm]{Remark}

\newtheorem{defn}{Definition}[section]

\title{When do Ten Points Lie on a Quadric Surface?}
\author{Will Traves}

\begin{document}
\maketitle

\begin{abstract}
A solution is provided to the Bruxelles Problem, a geometric decision problem originally posed 
in 1825, that asks for a synthetic construction to determine when ten points in 3-space lie on a quadric surface, a surface given by the vanishing of a degree-2 polynomial. The solution constructs four new points that are coplanar precisely when the ten original points lie on a quadric surface. The solution uses only lines constructed through two known points, planes constructed through three known points, and intersections of these objects. The tools involved include an extension of the Area Principle to three-dimensional space, bracket polynomials and the Grassmann-Cayley algebra, and von Staudt's results on geometric arithmetic. Many special cases are treated directly, leading to the generic case, where three pairs of the points generate skew lines and the remaining four points are in general position. A key step in the generic case involves finding a nice basis for the quadrics that pass through six of the ten points, which uses insights derived from {\it Macaulay2}, a computational algebra package not available in the nineteenth century. 
\end{abstract}

Dedicated to Bernd Sturmfels in recognition of his mathematical contributions and generous mentorship.

\section{Introduction} \label{ss:intro}

In 1825, l'Acad\'{e}mie de Bruxelles posed a geometric decision problem, asking for a synthetic construction that determines when ten points in 3-space lie on a quadric surface, a surface given by the vanishing of a degree-2 polynomial. Sturmfels highlighted this unsolved problem in his terrific book, Algorithms in Invariant Theory \cite[Exercise 3 on page 110]{S-AIT}. In this paper, we solve the Bruxelles Problem, constructing four new points that are coplanar precisely when the 10 original points lie on a quadric surface. 

The solution uses only lines constructed through two known points, planes constructed through 3 known points, and intersections of these objects. These are examples of the meet and join operations of the Grassmann-Cayley algebra, which we describe in Section \ref{ss:GC}. For our purposes, the join of two linear spaces in general position is just their span, and the meet of two linear spaces in general position is just their intersection. So the join of two distinct points $a$ and $b$ is the line $ab$ and the join of three non-collinear points $a$, $b$, and $c$ is the plane $abc$ containing all the points. 

The combined work of Pascal, Braikenridge, and Maclaurin \cite[Section 3.8]{CG} solved the much simpler version of the Bruxelles Problem for conic curves. 

\begin{figure}[h!t]
\begin{tikzpicture}

\tikzset{
every line/.append style ={
black,fill=black,line width=1pt}
} 

    \draw (0,0) ellipse [x radius=2, y radius=1];
    \draw [black,fill=black] (0,1) circle (.05cm);
    \draw [black,fill=black] (-1,0.86) circle (.05cm);
    \draw [black,fill=black] (1.2,0.79) circle (.05cm);
    \draw [black,fill=black] (1.3,-0.76) circle (.05cm);
    \draw [black,fill=black] (0.2,-0.99) circle (.05cm);
    \draw [black,fill=black] (-1.2,-0.8) circle (.05cm);
    \draw  (0,1)--(1.3,-0.76);
    \draw  (0,1)--(-1.2,-0.8);
    \draw  (1.3,-0.76)--(-1,0.86);
    \draw  (0.2,-0.99)--(-1,0.86);
    \draw  (0.2,-0.99)--(1.2,0.79);
    \draw  (-1.2,-0.8)--(1.2,0.79);
    \draw [black,fill=black] (-0.56,0.17) circle (.05cm);
    \draw [black,fill=black] (0.12,0.075) circle (.05cm);
    \draw [black,fill=black] (0.75,0.01) circle (.05cm);
    \draw [dashed] (-1.4,0.27)--(1.4,-0.08);
    \node  at (0,1.25) {$2$};
    \node  at (-1,1.11) {$0$};
    \node  at (1.2,1.04) {$4$};
    \node  at (1.3,-1.01) {$1$};
    \node  at (0.2,-1.24) {$5$};
    \node  at (-1.2,-1.05) {$3$};
    \end{tikzpicture}   
\caption{Six points lie on a conic when the points $05\cap 23$, $01\cap 34$, and $45\cap 12$ are collinear.}
\label{fig:pascal}
\end{figure}
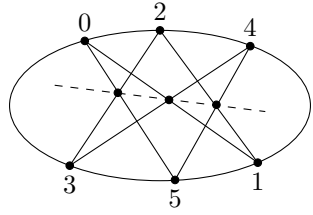

\begin{thm}[Pascal, Braikenridge-Maclaurin] \label{thm:PBM} Six points labeled $0$, $1$, $2$, $3$, $4$ and $5$ in the plane lie on a conic curve precisely when the three points $05\cap 23$, $01\cap 34$, and $45\cap 12$ are collinear, as illustrated in Figure \ref{fig:pascal}. \end{thm}

When the six points are the vertices of a regular hexagon in cyclic order, the intersecting lines in the theorem are parallel and meet on the line at infinity. It is thus natural to interpret problems in synthetic geometry as set in projective space. Projective space was already in use when the Bruxelles Problem was posed but it was not until 1847 that it was given a full modern treatment, in von Staudt's series of books on the Geometry of Position (1847-1860)\footnote{Geometrie der Lage (1847) and three volumes of Beitr\"{a}ge zur Geometrie der Lage (1856, 1857, 1860).}. Nevertheless, the modern version of the Bruxelles Problem sets the 10 points in projective 3-space, $\mathbb{P}^3$.

Many geometers of the nineteenth century attempted to solve the Bruxelles Problem, including Chasles, Hesse, Serret, and Steiner \cite[footnotes on first page]{TY}.  In his Ausdehnungslehre of 1862 \cite[Section 324]{HG}, Grassmann gave an elegant solution to the Bruxelles Problem when 3 skew lines, $L_0$, $L_1$, and $L_2$, each contain 3 of the ten points. In this case, the tenth point $x$ lies on the unique quadric surface containing the lines if the plane $xL_0$ meets the line $L_1$ in a point that lies on the plane $xL_2$. He uses the compact notation $xL_1L_2L_3x=0$ to describe this situation. In Lemma \ref{lem:lines}, we extend Grassmann's result to the setting where two lines each contain 3 of the ten points. Ciliberto and Miranda \cite{CM} consider many special cases of the Bruxelles Problem, characterizing when 10 points lie on a quadric using rational normal curves.  

The Bruxelles Problem is computationally feasible: the ten points $[x_i:y_i:z_i:w_i]$ ($i \in \{0,\ldots,9\}$) lie on a quadric surface precisely when the double Veronese images of the ten points lie on a hypersurface; that is, when the determinant of the $10 \times 10$ matrix  
$$N = \left[ \begin{array}{llllllllll} x_0^2 & x_0y_0 & x_0z_0 & x_0w_0 & y_0^2 & y_0z_0 & y_0w_0 & z_0^2 & z_0w_0 & w_0^2  \\ 
x_1^2 & x_1y_1 & x_1z_1 & x_1w_1 & y_1^2 & y_1z_1 & y_1w_1 & z_1^2 & z_1w_1 & w_1^2  \\
& & & & & \vdots & & & & \\
x_9^2 & x_9y_9 & x_9z_9 & x_9w_9 & y_9^2 & y_9z_9 & y_9w_9 & z_9^2 & z_9w_9 & w_9^2
\end{array} \right]$$
vanishes. 

The property that ten points lie on a quadric is invariant under linear change of coordinates so the Fundamental Theorem of Invariant Theory implies that this determinant can be written as a polynomial in the brackets $[abcd]$, where $[abcd]$ stands for the determinant of the $4\times 4$ matrix whose columns consist of the homogeneous coordinates of four points $a$, $b$, $c$, and $d$. Turnbull and Young \cite{TY} used invariant theoretic methods to construct a 240-term bracket polynomial that vanishes precisely when the 10 points lie on a quadric surface. White \cite{White90} showed that their polynomial can be straightened to a shorter polynomial with 138 terms. Apel and Richter-Gebert    \cite{SA, JRGSA} showed that the vanishing of bracket polynomials in $\mathbb{P}^2$ can be interpreted in terms of oriented distance ratios and established the Area Principle, a method for translating the vanishing condition for generic collections of points in terms of the Grassmann-Cayley algebra. One key step in our solution, described in Lemma \ref{lem:L.P} and Remark \ref{rmk:L.P}, extends their techniques to 3-dimensional space. 

Sturmfels and Whiteley \cite{SW} proved that every bracket polynomial $P$ admits a bracket monomial, a product of brackets $T$, so that $TP$ can be Cayley factored: $TP$ can be expressed using a combination of join and meet operations. Though $P$ may not factor directly, a factorization exists for almost all point configurations, those for which the bracket monomial $T$ does not vanish. Such Cayley expressions using just the meet and join operations immediately provide a synthetic construction to check whether the expression vanishes. The main step in our construction, described in Section \ref{ss:generic}, applies when a degree-4 bracket monomial does not vanish. In this generic case, the construction produces 4 test points in $\mathbb{P}^3$ that are coplanar precisely when the 10 points lie on a quadric surface. 

We give a high level overview of our construction. The pullback of the double Veronese map $\nu_2:\mathbb{P}^3 \rightarrow \mathbb{P}^9$ induces an isomorphism $\nu_2^*: (\text{Sym}^2(\mathbb{C}^{4}))^* \rightarrow \text{Sym}^2(\mathbb{C}^4{*})$, where 
linear functions on $\mathbb{P}^9$ are sent to quadratic functions on $\mathbb{P}^3$. We decompose $ (\text{Sym}^2(\mathbb{C}^{4}))^* \cong \mathbb{C}^{10*} = S \oplus T$, where $S$ is the vector subspace of linear functions vanishing on points 
$0,\ldots,5$ and $T$ is its complement. Then we consider the projection $\pi: \mathbb{C}^{10*} \rightarrow S$ that kills $T$. 
The injection map splits: the map $i:S \hookrightarrow S\oplus T = \mathbb{C}^{10*}$ satisfies $\pi\circ 
i = \text{id}_S$. 
In the nongeneric case, $\text{dim}(T)<6$, so the Veronese images of the first six points fail to impose linearly independent conditions and so there is a degree-2 function vanishing on the 10 points; that is, the 10 points lie 
on a quadric surface. Determining if we are in the nongeneric case using only synthetic tools is the content of Section \ref{s:10pts}. In essence, cases where the genericity assumptions fail are so constrained that we can answer the Bruxelles Problem directly. In the generic case, $\text{dim}(T)=6$ 
and we find a nice basis for $S$, consisting of linear functions corresponding to reducible quadratic functions on $\mathbb{P}^3$. Finding a basis depends on insight developed working with the computer algebra system \emph{Macaulay2} \cite{M2}. In 
this case, the Veronese images of the last four points fail to impose linearly independent conditions precisely when a $4\times 4$ matrix $M$ is singular. The vanishing of $\text{det}(M)$ implies a nonzero linear function $s \in S$ which induces a quadratic function $\nu_2^*(i(s))$ vanishing on all 10 points. 
Constructing the matrix $M$ and 
determining if its determinant is zero using synthetic methods is the main content of Section \ref{ss:generic}.  
The argument uses that three lines, each generated by pairs of our points, are mutually skew and that the remaining four points are in general position (not lying in a hyperplane); we reduce to this case in Section \ref{s:10pts}. 

Each column of the $4 \times 4$  matrix $M$ corresponds to a point in $\mathbb{P}^3$: the 10 points lie on a quadric surface when the four test points $P_0$, $P_1$, $P_2$, and $P_3$ are coplanar. In Section \ref{ss:test} we show how the coordinates of each test point can be obtained as the product of a quotient of brackets if we work in an appropriate chart for each test point. Using the extension of the Area Principle to 3-dimensions, we identify the quotients with points on a line and construct a point on the line corresponding to their product using constructions originally due to von Staudt (and described in Section \ref{ss:prod}). The product points for three coordinates enable us to reconstruct the image $Q_i=\tau(P_i)$ of the test point $P_i$  under an automorphism $\tau$ of $\mathbb{P}^3$, and the coplanarity of these four image points $Q_0$, $Q_1$, $Q_2$, and $Q_3$ tests whether the original 10 points lie on a quadric surface. 

The article is organized into three additional sections. Section \ref{s:incidencebg} gives important background information. Our solution depends on defining local coordinates on six lines given by extending the edges in a tetrahedron in $\mathbb{P}^3$. Local coordinates and cross-ratios are reviewed in Section \ref{ss:cr}. We can recover points in $\mathbb{P}^3$ from their projections onto (some of) these six lines, a process described in Section \ref{ss:recovering}.  Section \ref{ss:generic} describes our synthetic construction in the generic case. Section \ref{s:10pts} explains how to handle the non-generic cases and reduce to the generic case. The arguments in this section involve computations in the Grassmann-Cayley algebra, which is described in section \ref{ss:GC}. We advertise some open problems and provide suggestions for further reading in the last section.

\section{Incidence in Projective Space}
\label{s:incidencebg}
Every algebraic geometer knows that we need to fix 3 points to kill automorphisms of $\mathbb{P}^1$ and that any $n+2$ points in linearly general position in $\mathbb{P}^n$ can be sent to any other $n+2$ points in linearly general position via a unique automorphism of $\mathbb{P}^n$. Both results play key roles in our solution to the Bruxelles Problem. 

\subsection{Cross ratios as local coordinates} \label{ss:cr}

We first recall some additional facts about the projective line. After a linear change of coordinates on the underlying two-dimensional space, we can parameterize the line so that three given points -- zero, unit, and infinity -- are located at  $p_0=[1:0]$, $p_1=[1:1]$, and $p_\infty=[0:1]$, respectively. The cross ratio $(a,b;c,d)$ of any four labeled points on $\mathbb{P}^1$ is the point $p=\phi(d)$, where $\phi$ is the unique linear automorphism of $\mathbb{P}^1$ that sends $a$ to $p_\infty$, $b$ to $p_0$, and $c$ to $p_1$, so that $(p_\infty,p_0;p_1,p)=p$. The 
point $p_x=[1:x]$ can be identified with the parameter value $x$ and this parameter can be recovered using determinants to compute the cross ratio \cite[Lemma 5.5]{RG}:
$$(p_\infty,p_0; p_1,p_x) = \frac{[p_\infty p_1][p_0 p_x]}{[p_\infty p_x][p_0 p_1]} = \frac{(-1)(x)}{(-1)(1)}=x ,$$ where the notation $[\cdot]$ stands for the determinant of the square matrix whose columns are the homogeneous coordinates of the points enclosed by the bracket. The computation of the cross ratio does not depend on the choice of representatives for the points: multiplying any of the points by a constant does not change the cross ratio. Furthermore, the value of the cross ratio is invariant under change of coordinates on the projective line since if $A \in GL_2$, $[(Ap_i)(Ap_j)] = \det(A)[p_ip_j]$.  When the four points $p_0$, $p_1$, $p_\infty$, and $p_x$ lie on an arbitrary line in $\mathbb{P}^2$, we can compute their cross ratio using any point $q$ off the line, $$(p_\infty,p_0; p_1,p_x) = \frac{[p_\infty p_1 q][p_0 p_x q]}{[p_\infty p_x q][p_0 p_1 q]}.$$
This is easy to check if we choose coordinates so that $p_0$, $p_\infty$, and $q$ are the standard basis vectors  for the vector space underlying $\mathbb{P}^2$. Similarly, 
when the points are collinear in $\mathbb{P}^3$ and $q$ and $r$ are two distinct points off the line, 
$$(p_\infty,p_0; p_1,p_x) = \frac{[p_\infty p_1 qr][p_0 p_x qr]}{[p_\infty p_x qr][p_0 p_1 qr]}.$$

\subsection{The product construction} \label{ss:prod}  
\begin{figure}[h!t]
\begin{center}
\begin{tikzpicture}

\tikzset{
every line/.append style ={
black,fill=black,line width=1pt}
} 
    \draw [black,fill=black] (0,0) circle (.05cm);
    \draw[black,fill=black] (2,0) circle (0.05cm);
    \draw[black,fill=black] (3,0) circle (0.05cm);

    \draw[black,fill=black] (5.35,0) circle (0.05cm);
    \draw[black,fill=black] (7.04,0) circle (0.05cm);
    \draw[black,fill=black] (14,0) circle (0.05cm);
    \draw (-0.5,0)--(14.5,0);
    \draw (0,-0.5)--(0,4.5);
    \draw[black,fill=black] (0,2.16) circle (0.05cm);
\draw[black,fill=black] (0,3) circle (0.05cm);
\draw[black,fill=black] (0,4) circle (0.05cm);
\draw[black,fill=black] (1,2) circle (0.05cm);
\draw[black,fill=black] (3.11,1.68) circle (0.05cm); 
\draw (0,4)--(5.34,0);
\draw (0,3)--(7.04,0);
\draw (0,3)--(3,0);
\draw (0,4)--(2,0);
\draw (0,2.16)--(14,0); 
\node at (1.15,2.15)
{$a$};
\node at (3.26,1.87)
{$b$};
\node at (14.2,-0.3)
{$p_{\infty}$};
\node at (7.2,-0.3)
{$p_y''=p_{xy}$};
\node at (0.3,-0.3)
{$p_0$};
\node at (2.2,-0.3)
{$p_1$};
\node at (3.2,-0.3)
{$p_y$};
\node at (5.55,-0.3)
{$p_x$};
\node at (-0.3,4)
{$p_1'$};
\node at (-0.3,3)
{$p_y'$};
\node at (-0.3,2.16)
{$p_\infty'$};
\node at (0.3,4.5)
{$L'$};
\node at (14.5,0.3)
{$L$};
\end{tikzpicture} 
\end{center}
\caption{The product construction.}
\label{fig:multx}
\end{figure}
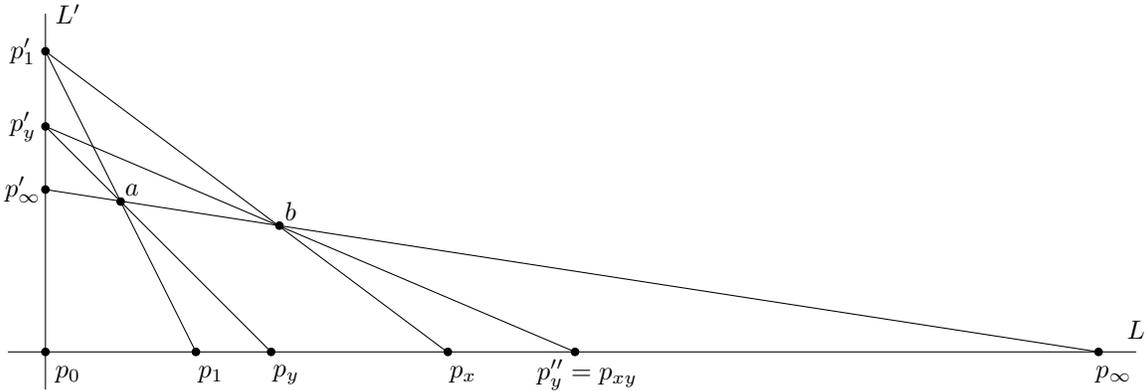

Consider a line $L$ in $\mathbb{P}^2$. After a change of basis on $\mathbb{P}^2$ we may assume that $L$ is the line spanned by $p_0=[1:0:0]$ and $p_\infty=[0:1:0]$ and the unit on $L$ is $p_1=[1:1:0]$. Given a point $p_x=[1:x:0]$ on $L$, the map $\gamma: L \rightarrow L$ fixing $p_0$ and $p_\infty$ and sending $p_1$ to $p_x$ corresponds to left multiplication by the $3\times 3$ diagonal matrix $A_x = \text{diag}(1,x,1)$ which sends $p_y=[1:y:0]$ to $p_{xy}=[1:xy:0]$. That is, scaling the unit by $x$, multiplies $p_y$ by $x$ as well. The mathematician von Staudt showed that the map $\gamma$ can be achieved by a composition of two maps, each given by projection from a point, as illustrated in Figure \ref{fig:multx}. We start by picking a point $a$ off $L$ and another line $L'$ through $p_0$. Define the points $p_0'=p_0$, $p_1'$, $p_y'$, and $p_\infty'=p_\infty$ on $L'$ to be the projections of $p_0$, $p_1$, $p_y$, and $p_\infty$ from $a$ onto $L'$. Now set $b$ equal to the intersection of the line joining $p_1'$ to $x$ with the line joining $a$ and $p_\infty$. Then define the points $p_0^{\prime\prime}=p_0$, $p_1^{\prime\prime}$, $p_y^{\prime\prime}$, and $p_\infty^{\prime\prime}=p_\infty$ on $L$ to be the projections of $p_0'$, $p_1'$, $p_y'$, and $p_\infty'$ from $b$ onto $L$. The projection of $L$ from $a$ onto $L'$ followed by the projection of $L'$ from $b$ onto $L$ fixes $p_0$ and $p_\infty$ and sends $p_1$ to $p_x$, so the composition equals the map $\gamma$; in particular, $p_{y}^{\prime\prime}=p_{xy}$ is constructed from $p_0$, $p_1$, $p_\infty$, $p_x$ and $p_y$ in this manner. Though we chose a convenient basis to interpret the  construction, the construction of the point corresponding to $xy$ from the points corresponding to $x$ and $y$ works for any line $L$ in $\mathbb{P}^2$, and even for lines in higher-dimensional projective spaces.  Note that this product construction only depends on incidence relations and the positions of the points $p_0$, $p_1$, and $p_\infty$ on $L$. 

Reversing the order of the construction in an appropriate way gives a construction of $p_{1/x}$ from $p_{x}$. First we draw a new line $L'$ through $p_0$. Then we pick a point $b$ off $L \cup L'$ and draw lines $p_{\infty}b$, $p_x b$ and $p_1b$. We draw the line $(p_xb\cap L')p_1$ and let $a$ denote the intersection of this line and $p_\infty b$. Then we draw $(p_1b\cap L')a$, which intersects $L$ at the point $p_{1/x}$. This construction works fine if $x=0$ or $x=\infty$: if $x=0$ it produces $p_\infty$ and if $x=\infty$ it produces $p_0$. 

\subsection{Recovering points from local coordinates} 
\label{ss:recovering} Now we focus on 3-dimensional projective space. 
After a change of coordinates on the underlying vector space, any five points in linearly general position in $\mathbb{P}^3$ (no four on a plane) are located at $E_0=[1:0:0:0]$, $E_1=[0:1:0:0]$, $E_2=[0:0:1:0]$, $E_3=[0:0:0:1]$ and $\mathbf{1}=[1:1:1:1]$.
If $v$ is a nonzero vector in four-dimensional space, write $\langle v \rangle$ for the point in $\mathbb{P}^3$ corresponding to the 1-dimensional subspace containing $v$. Then the point $E_i = \langle e_i \rangle$ corresponds to the standard basis vector $e_i$ and $\mathbf{1} = \langle e_0+e_1+e_2+e_3 \rangle$. When $i<j$ we take $E_i$ and $E_j$ as the points zero and infinity on $E_iE_j$. Define $\mathbf{1}_{ij}$ to be the point $\langle e_i+e_j \rangle$ on the line $E_iE_j$; $\mathbf{1}_{ij}$ serves as the unit on $E_iE_j$. 

The unit points on the lines $E_iE_j$ are consistently determined by $\mathbf{1}$ through an incidence construction.  If $\{i,j,k,\ell\} = \{0,1,2,3\}$, define the projection $\pi_{ij}: \mathbb{P}^3 \dashrightarrow E_iE_j$ sending $P$ to the intersection of the line $E_iE_j$ with the plane $E_k E_\ell P$ containing $E_k$, $E_\ell$ and $P$. 
The projection morphism $\pi_{ij}$ is  defined off the line $E_kE_\ell$.  
The point $P=[x_0:x_1:x_2:x_3] = \langle x_0e_0+x_1e_1+x_2e_2+x_3e_3 \rangle$ is sent by the projection $\pi_{ij}$ to $\pi_{ij}(P) = \langle x_ie_i + x_je_j\rangle$. In particular, $\pi_{ij}(\mathbf{1}) = \mathbf{1}_{ij}$. If $i<j$, then $\pi_{ij}(P)$ is the point on $E_iE_j$ corresponding to the parameter value $x_j/x_i$. In this way, the locations of the five points $E_0$, $E_1$, $E_2$, $E_3$, and $\mathbf{1}$ in $\mathbb{P}^3$ determine units and local parameters on all six lines $E_iE_j$.

For $0\leq i \leq 3$, let $U_i$ stand for the chart $\{[x_0:x_1:x_2:x_3] \in \mathbb{P}^3:\, x_i \neq 0\} \cong K^3$. The four charts $U_i$  cover $\mathbb{P}^3$. We say that the chart $U_i$ is spanned by the lines $E_{i}E_{j}$ with $j\neq i$ and that the line $E_iE_j$ is part of a spanning set for $U_i$. We can use projections to determine if a point lies in a chart: the point $P$ lies in $U_i$ if and only if the projections $\pi_{ij}(P)$ exist and are not equal to the infinity point for all three lines $E_iE_j$ spanning $U_i$. Moreover, given the projections of a point $P$ onto the lines spanning a chart $U_i$, we can recover $P$ by intersecting three transverse planes: if $\{i,j,k,\ell\}=\{0,1,2,3\}$ and  $P$ lies in $U_i$ then 
$$P = \pi_{ij}(P)E_kE_\ell \cap
\pi_{ik}(P)E_jE_\ell \cap
\pi_{i\ell}(P)E_jE_k.$$

\subsection{The Grassmann-Cayley Algebra} \label{ss:GC}

We give a brief introduction to the meet and join operations in the Grassmann-Cayley algebra; for details, see Sturmfels \cite[Section 3.3]{S-AIT} Given a $d$-dimensional vector space $V$, the meet and join operations define product structures on the exterior algebra $\oplus_{k=0}^d \bigwedge^k V$. 
We say that an element is a $k$-extensor if it is the wedge product of $k$ vectors. The two operations are linear, so it suffices to define them for products of extensors. We warn the reader that the usual notation for meets ($\wedge$) and joins ($\vee$)  conflicts with the wedge product: we write $A = a_1\vee a_2 \vee \cdots a_k$ for the wedge product of $k$ vectors $a_i$ and abbreviate this expression as $A =a_1a_2\cdots a_k$. A nonzero extensor $A=a_1a_2\cdots a_k$ 
involves vectors $a_1,\ldots,a_k$ that span a $k$-dimensional vector space $\bar{A} = \{v \in V: a_1a_2\cdots a_kv=0\} $. Conversely, any $k$-dimensional vector space $\bar{A}$ determines (up to scalar multiple) an extensor $A$ given by the wedge product of a basis for $\bar{A}$.  
If $A=a_1a_2\cdots a_j$ is a $j$-extensor and $B=b_1b_2\cdots b_k$ is a $k$-extensor, then we their join is a $j+k$-extensor $A \vee B = a_1a_2\cdots a_jb_1b_2\cdots b_k$. If the join $A \vee B$ is nonzero, then the vectors $a_1,\ldots,a_j,b_1,\ldots,b_k$ are linearly independent and $\overline{A\vee B}$ equals $\overline{A}+\overline{B}$, the span of $\overline{A}$ and $\overline{B}$. Fix a basis $e_1,\ldots,e_d$ for the vector space $V$ and identify $\wedge^d V$ with the 1-dimensional vector space spanned by this basis, sending $e_1\cdots e_d$ to $1$. Then any $d$-extensor $a_1\cdots a_d$ can be expressed as the determinant $[a_1\cdots a_d]$.   

The meet $A \wedge B$ of $A=a_1a_2\cdots a_j$ and $B=b_1b_2\cdots b_k$ defined to be zero if $j+k < d$ and 
$$ A \wedge B = \sum_\sigma \text{sign}(\sigma) [a_{\sigma(1)}\cdots a_{\sigma(d-k)}b_1\cdots b_k]a_{\sigma(d-k+1)}\cdots a_{\sigma(j)},$$
where the sum is over all $(d-k,j-(d-k))$-split shuffles, permutations of $\{1,\ldots, j\}$ with $\sigma(1)<\cdots<\sigma(d-k)$ and $\sigma(d-k+1)<\ldots<\sigma_d$. The meet of two extensors is an extensor and we have $A\wedge B$ is nonzero if and only if $\overline{A}+\overline{B}=V$. In that case, the meet corresponds to the intersection of subspaces, $\overline{A \wedge B} = \overline{A} \cap \overline{B}$. 

If $p$, $q$ and $r$ are points in a projective space, then the join $pq$ corresponds to the line through $p$ and $q$ ($pq=0$ if $p=q$). The join $pqr$ corresponds to the plane containing $p$, $q$ and $r$ ($pqr=0$ if the points lie on a line). In  $\mathbb{P}^2$ the meet can be used to give an explicit description of the point of intersection of two lines, $$ pq \wedge rs = [prs]q-[qrs]p$$ and in $\mathbb{P}^3$ we can find explicit expressions for the intersection of a line with a plane, 
$$pq\wedge rst =[prst]q - [qrst]p,$$
and the intersection of two planes, 
$$ pqr\wedge stu  =  [pstu]qr - [qstu]pr + [rstu]pq.$$

\section{The generic case} \label{ss:generic}

In this section, we describe an algorithm to decide whether the 10 points $0$, $1$, $\ldots$, $9$ lie on a quadric surface. We assume that the 10 points satisfy four genericity conditions:
\begin{enumerate}
\item all ten points are distinct; 
\item no four of the ten points are collinear; 
\item the lines $01$, $23$, and $45$ are mutually skew; 
\item the four points $6$, $7$, $8$, and $9$ are in general position (not contained in a plane). 
\end{enumerate}

Note that the ten points satisfy the genericity conditions if and only if the bracket polynomial $[0123][0145][2345][6789]$ is nonzero. 
In Section \ref{s:10pts} we reduce the general problem to this crucial case. For now, we note that if conditions (1) or (2) fail then the 10 points fail to impose independent conditions on quadrics and so they lie on a quadric surface. 

\subsection{A family of quadrics}

We focus on the quadric surfaces that pass through the six points $0$, $1$, $2$, $3$, $4$, and $5$. We claim that the assumption that $01$, $23$ and $45$ are skew lines implies that the six points impose independent conditions on quadrics, so there is a $\mathbb{P}^3$ of quadrics through all six points. Since the points $0$, $1$, $2$ and $3$ are in general position, we can apply a projective linear transformation to the points to put them at $[1:0:0:0]$, $[0:1:0:0]$, $[0:0:1:0]$, and $[0:0:0:1]$. Suppose that point $4$ is at $[x_4:y_4:z_4:w_4]$ and point $5$ is at $[x_5:y_5:z_5:w_5]$. Representing a general quadric with the equation $$ax^2 + bxy + cxz + dxw + ey^2 + fyz + gyw + hz^2 + izw + jw^2 =0,$$ we see that the quadrics passing through points $0$, $1$, $2$, and $3$ have $a=e=h=j=0$ and the points $4$ and $5$ impose dependent conditions when 
$$
\text{Rank}\left( \begin{bmatrix} x_4y_4 & x_4z_4 & x_4w_4 &  y_4z_4 & y_5w_5 & z_5w_5  \\
x_5y_5 & x_5z_5 & x_5w_5 &  y_5z_5 & y_5w_5 & z_5w_5  \end{bmatrix}\right) \leq 1.$$ This occurs when all the $2 \times 2$ minors of this matrix vanish. We used \emph{Macaulay2} \cite{M2} to find the primary decomposition of the ideal generated by the minors. The primary components correspond to one of the following conditions: points $4$ and $5$ are equal; point $4$ or $5$ equals one of points $0$, $1$, $2$, or $3$; or the points $4$ and $5$ are collinear with two of the four points $0$, $1$, $2$, and $3$. All of these situations violate our assumptions for this section, so we conclude that the six points $0$, $1$, $2$, $3$, $4$, and $5$ impose independent conditions on quadrics and the set of quadrics passing through the six points forms a $\mathbb{P}^3$. 

We now find a basis for the quadrics passing through the first six points. 
Write $[abcX]$ for the determinant of the $4 \times 4$ matrix whose columns are the points $a$, $b$, $c$, and the variable point $X = [x:y:z:w]$ so that $[abcX]=0$ is the equation of the plane through points $a$, $b$, and $c$. We pick a set of four reducible quadrics through the points $0, \ldots, 5$, 
$$ \mathcal{S} = \{ [015X][234X], [012X][345X], [024X][135X], [045X][123X] \}$$
and investigate when these span the space of quadrics through all six points. It is easier to check the equivalent condition for these four polynomials to be linearly independent. Working in the ring $\mathbb{Q}[x_4,y_4,z_4,w_4,x_5,y_5,z_5,w_5][x,y,z,w]$ we compute the $4 \times 10$ matrix whose entries in each column are the coefficients of the degree-2 monomials in $x$, $y$, $z$, and $w$ appearing in the polynomials in $\mathcal{S}$. The polynomials in $\mathcal{S}$ are linearly dependent when all the $4 \times 4$ minors of this matrix vanish. Using \emph{Macaulay2} \cite{M2}, we find that the ideal generated by these minors is a primary ideal and its minimal prime ideal is generated by 
$$ Q = -x_5y_4z_5w_4 + x_4y_5z_5w_4
+ x_5y_4z_4w_5 - x_4y_4z_5w_5. $$

We give an incidence interpretation of the condition $Q=0$. 

\begin{figure}[h!t]
\begin{center}
    \begin{tikzpicture}[x=0.75in,y=0.75in]

\tikzset{
every line/.append style ={
black,fill=black,line width=1pt}
}   
    \draw [black,fill=black] (0,0) circle (.05cm);
    \draw [black,fill=black] (1.2,1) circle (.05cm);
    \draw [black,fill=black] (2,0) circle (.05cm);
    \draw [black,fill=black] (1.6,0.51) circle (.05cm);
    \draw [black,fill=black] (0.5,0.42) circle (.05cm);
    \draw (0,0)--(1.6,0.51);
    \draw (2,0)--(0.5,0.42);
    \draw [black,fill=black] (0.94,0.3) circle (.05cm);
    \draw [black,fill=black] (0.84,0) circle (.05cm);
    \draw  (0.84,0)--(1.2,1);
    \draw  (0,0)--(1.2,1);
    \draw  (0,0)--(2,0);
    \draw  (2,0)--(1.2,1);
    \node at (0,-0.2) {$2$};
    \node at (2,-0.2) {$3$};
    \node  at (1.2,1.2) {$1$};
    \node at (0.84,-0.2) {$p=23 \cap 015$};
    \node at (2.2,0.545) {$q=13\cap 024$};
    \node at (-0.1,0.46) {$r=12 \cap 345$};
\end{tikzpicture} 
\end{center}
\vspace{-0.1in}
\caption{An incidence interpretation of $Q=0$.} 
\label{fig:Ceva-like}
\end{figure}
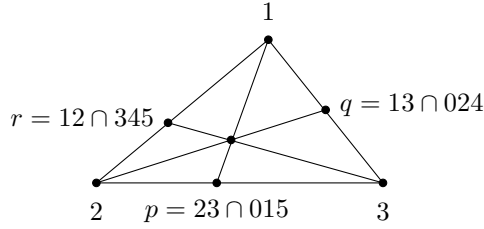

\begin{lem} \label{lem:Q} The expression $Q$ is equal to the evaluation of the bracket binomial 
$$ [0125][0234][1345] - [0124][2345][0135]$$ 
and both vanish precisely when an incidence condition reminiscent of Ceva's Theorem holds: 
if we set $p = 23 \cap 015$ equal to the intersection of the line $23$ and the plane $015$, $q=13\cap 024$ and $r = 12 \cap 345$ then the condition $Q=0$ requires that the lines $1p$, $2q$, and $3r$ meet in a common point, as shown in Figure \ref{fig:Ceva-like}. 
\end{lem}

\begin{proof}
A direct computation shows that the bracket polynomial equals $Q$. Taking the meet of $23$ with $015$ gives an explicit description of $p=[2015]3-[3015]2 = [0125]3 - [0135]2$ as a point on the line $23$. Similarly, $q = -[0124]3-[0234]1$ and $r=[1345]2-[2345]1$. The lines $1p$ and $2q$ should meet because they are in the plane $123$ but this can be difficult to determine directly in the Grassmann-Cayley algebra. Instead, we cone over the point $0$ (which lies off the plane $123$ since $01$ and $23$ are skew) and intersect the planes $01p$ and $02q$, giving a line containing $0$ and the intersection point of $1p$ and $2q$. Taking the meet of $01p$ and $02q$ gives $-[102q]0p +[p02q]01$. The line corresponding to this element meets the line $3r$ only if the line $3r$ contains the intersection point of $1p$ and $2q$. Taking the meet of the line with $3r$ gives a bracket polynomial $-[102q][0p3r]+[p02q][013r]$ which must vanish precisely when $1p$, $2q$, and $3r$ are coincident. Substituting for $p$, $q$, and $r$, we obtain 
$$ -[102q][0p3r]+[p02q][013r] =[0123]^2 Q. $$ 
Since $01$ and $23$ are skew, $[0123]$ is nonzero, so $Q=0$ characterizes when $1p$, $2q$, and $3r$ are coincident.
\end{proof}

It may happen that the expression $Q$ is zero and the four quadrics in $\mathcal{S}$ are linearly dependent. We get around this difficulty by permuting the point labels. Given a permutation $\sigma$ of the set $\{0,1,2,3,4,5\}$, we set
$ \sigma([abcd]) = [\sigma(a)\sigma(b)\sigma(c)\sigma(d)]$ and extend this action to a ring homomorphism on the bracket algebra (the polynomial ring on the brackets modulo the ideal generated by the relations among the brackets, the Grassmann-Pl\"{u}cker relations, which can be obtained in \emph{Macaulay2} \cite{M2} using the \emph{Grassmannian} command). 
Now we form the ideal $I = (\sigma(Q): \sigma \in \mathbb{S}_6)$, generated by all the permutations applied to the expression $Q$, thought of as an element of the bracket algebra. This ideal contains the polynomial $[0123][0145][2345]$, which vanishes only when two of the three lines $01$, $23$ and $45$ meet. Since the lines are skew, $\sigma(Q)\neq 0$ for some $\sigma$ and we can identify the permutation $\sigma$ using the incidence interpretation given in Lemma \ref{lem:Q}. Relabeling each of the six points $\sigma(i)$ with the label $i$, we find that $Q \neq 0$ and the polynomials in $\mathcal{S}$ are a basis for the space of all quadrics that pass through the six points. We note that after relabeling, the lines $01$, $23$, and $45$ may no longer be skew.

Searching for a quadric through all ten points requires us to find a polynomial of the form  
$$\alpha [015X][234X]+\beta [012X][345X]+\gamma [024X][135X]+\delta [045X][123X]$$
that vanishes at $X = 6$, $7$, $8$ and $9$. So we are led to solve a system 
$$ M\mathbf{v} = \begin{bmatrix} 
[0156][2346] &  [0126][3456] &  [0246][1356] &  [0456][1236] \\
[0157][2347] &  [0127][3457] &  [0247][1357] &  [0457][1237] \\
[0158][2348] &  [0128][3458] &  [0248][1358] &  [0458][1238] \\
[0159][2349] &  [0129][3459] &  [0249][1359] &  [0459][1239]                 
\end{bmatrix} \begin{bmatrix} \alpha \\ \beta \\ \gamma \\ \delta \end{bmatrix} = \begin{bmatrix} 0 \\ 0 \\ 0 \\ 0 \end{bmatrix}.$$

We can assume that none of the columns of $M$ are zero; otherwise, the points $6$, $7$, $8$, and $9$ lie on a union of two planes, each passing through 3 of the remaining points, so the 10 points lie on a union of two planes. There is a nice connection between the two matrices $M$ and $N$ and the polynomial $Q$:  $\det(M)=Q\det(N)$. It is not clear whether this is a coincidence or an indication of some more general result. 

The system $M\mathbf{v}=\mathbf{0}$ has a nonzero solution precisely when $\det M = 0$, that is, when the four test points $P_1$, $P_2$, $P_3$, and $P_4$, whose coordinates are given by the columns of $M$, are coplanar. Rather than constructing these four points, we construct their images under an invertible projective transformation using meet and join operations. 

\subsection{Constructing the test points}
\label{ss:test}

We use the data in the columns of $M$ to construct a point in the 3-space containing the original 10 points. The four points $6$, $7$, $8$, and $9$ span a $\mathbb{P}^3$ and the fifth point $u$ determines the global unit in the $\mathbb{P}^3$ containing the
10 original points. Choose vectors so that points $6$, $7$, $8$, $9$ and $u$ equal $\langle v_6\rangle$, $\langle v_7\rangle$, $\langle v_8\rangle$, and $\langle v_9\rangle$, and $\langle v_6+v_7+v_8+v_9\rangle$, respectively. The $\mathbb{P}^3$ containing the ten 
original points is not the same as the $\mathbb{P}^3$ spanned by the columns of $M$, though they are isomorphic.    We define an isomorphism $\tau$ between the two
$\mathbb{P}^3$'s, sending $E_0$, $E_1$, $E_2$, $E_3$, and $\mathbf{1}$
in the $\mathbb{P}^3$ generated by the 
columns of $M$ to the points $6$, $7$, $8$, $9$, and $u$, respectively.  

As an example of the method, consider $P_1$, the first column of $M$. The point $P_1$ lies in some chart of $\mathbb{P}^3$; for illustrative purposes, we suppose that $P_1$ lies in $U_0$. Dividing all the coordinates of $P_1$ by the coordinate in the $0^\text{th}$ position, we get an affine representative for $P_1$, 
$$ \left[1: \frac{[0157]}{[0156]}\frac{[2347]}{[2346]}: \frac{[0158][2348]}{[0156][2346]}:\frac{[0159][2349]}{[0156][2346]} \right].
$$
Now we focus on the second coordinate of this representation of $P_1$. Viewing the coordinate as the product of two fractions, we are inspired to construct a point on a line whose local coordinate is a given ratio of brackets sharing three common entries. 

\begin{lem} \label{lem:L.P} If $0 \leq a < b < c \leq 5$ and $6 \leq d < e \leq 9$ then the point $p = de \wedge abc = [dabc]e - [eabc]d = [abce]d-[abcd]e$ corresponds to the local parameter value $-[abcd]/[abce]$.   
\end{lem}

\begin{proof} The expression for $p$ in terms of $d$ and $e$ follows directly from the definition of the meet operation. The local parameter value can be recovered using the cross product. The point $d = \langle v_d \rangle$ serves as zero and the point $e = \langle v_e \rangle$ serves as infinity on the line $de$. The unit on $de$ is $de\wedge fgu$, where $\{d,e,f,g\} = \{6,7,8,9\}$ and $u=\langle v_6+v_7+v_8+v_9 \rangle = \langle v_d+v_e+v_f+v_g \rangle$. 
Computing, $v_dv_e\wedge v_fv_gv_u = [v_dv_fv_gv_u]v_e - [v_ev_fv_gv_u]v_d = [v_dv_fv_gv_e]v_e-[v_ev_fv_gv_d]v_d$, so the unit $\mathbf{1}_{de}$ on $de$ is $\langle [v_dv_ev_fv_g](v_e+v_d) \rangle = \langle v_d+v_e\rangle$. The cross ratio $(\mathbf{\infty}_{de},\mathbf{0}_{de}; \mathbf{1}_{de}, p)$ equals the local coordinate of the point $p$ on $de$: 
$$ \begin{array}{lll} \vspace{0.1in}(\mathbf{\infty}_{de},\mathbf{0}_{de}; \mathbf{1}_{de}, p) &= &
\frac{[e(d+e)fg][dpfg]}{[epfg][d(d+e)fg]} \\ \vspace{0.1in} & = & \frac{[edfg][dpfg]}{[epfg][defg]} \\ \vspace{0.1in}
& = & \frac{[edfg][d([abce]d-[abcd]e)fg]}{[e([abce]d-[abcd]e)fg][defg]}
\\ \vspace{0.1in}
& = & - \frac{[edfg][abcd][defg]}{[abce][edfg][defg]}
\\ & = & - \frac{[abcd]}{[abce]}.\end{array}$$
\end{proof} 

\begin{rmk} \label{rmk:L.P} We can also find a point $q$ on $de$ with local parameter value $-[abce]/[abcd]$. We just apply the construction of $p_{1/x}$ from $p_x$ given in section \ref{ss:prod} to the point $p$ in Lemma \ref{lem:L.P}.
\end{rmk}

Using this remark, let  $p_x = 67 \wedge 015$ be the point $p_x$ on the line $67$ with local coordinate $x = -[0156]/[0157]$ and construct $q_x = p_{1/x}$. Also let $p_y=67 \wedge 345$ be the point with local coordinate $y=-[2346]/[2347]$ and construct $q_y=p_{1/y}$. We use the product construction from section \ref{ss:prod} on the line $67$ to produce a point $Q_{67}$ corresponding to the local parameter $(-\frac{1}{x})(-\frac{1}{y}) = [0157][2347]/[0156][2346]$  on the line $67$.  Similarly, we construct the point $Q_{68}$ corresponding to local parameter $[0158][2348]/[0156][2346]$ on the line $68$. Finally, we construct the point $Q_{69}$ corresponding to local parameter $[0159][2349]/[0156][2346]$ on the line $69$. Now we construct the point 
$$Q_1 = Q_{67}89 \cap Q_{68}79 \cap Q_{69}78.$$ This is the image of $P_1$ under the isomorphism $\tau$ of $\mathbb{P}^3$.

The same construction can be applied to the other columns of $M$ to construct the points $Q_2$, $Q_3$, and $Q_4$ that are the images under $\tau$ of $P_2$, $P_3$, and $P_4$.  The four columns of $M$ are linearly dependent precisely when the four points $Q_1$, $Q_2$, $Q_3$, and $Q_4$ are coplanar. It follows that the 10 original points lie on a quadric surface precisely when the four constructed points $Q_1$, $Q_2$, $Q_3$, and $Q_4$ are coplanar.

\section{Geometric Reductions} 
\label{s:10pts}

In this section we handle all cases of the general Bruxelles Problem that violate one of our four assumptions in Section \ref{ss:generic}. Two of these are easy. If two of the ten points agree then two of the rows in the constraint matrix $N$ agree and so $\text{det}(N)=0$ and the points lie on a quadric. If four of the ten points lie on a line, then the constraints imposed by three of the points force the fourth point constraint to be satisfied, since any quadric that contains 3 collinear points must contain the entire line. 

\begin{defn}[Quadric Decidable (QD)]
\label{def:QD} Say that the 10 points are quadric-decidable (QD) if the points are in special position so that: 
\begin{enumerate}
\item We have a synthetic construction to determine if all 10 points lie on a quadric surface
\item We also have a synthetic construction that checks whether the points lie in our special position. 
\end{enumerate} 
\end{defn}

Our goal in this section is to show that either the 10 points are QD or there are three pairs of our points that generate mutually skew lines and the remaining four points are in general position. We begin by finding find three pairs of our points that generate mutually skew lines.  

\subsection{Finding three skew lines}
\label{ss:geomred}

We start by showing that there are at least three skew lines among the lines generated by the ten points. If every pair of lines $ij$ and $k\ell$ meet, then the line $01$ meets all the lines $23$, $24$, \ldots, $29$ so that the points $3$, $4$, \ldots, $9$ all lie on the plane $012$ (and hence on a quadric surface). So we can restrict attention to the case where there are at least two skew lines, $01$ and $23$, say. We can assume that at least one of the remaining points $4,5,\ldots,9$ lies off both lines $01$ and $23$, otherwise, the ten points lie on a quadric surface containing the lines $01$ and $23$. Relabel so that $4$ lies off both $01$ and $23$.   
Now we claim that one of the lines $45$, $\ldots$, $49$ is skew to both $01$ and $23$, otherwise all the points $5,6,\ldots,9$ must lie on the degenerate quadric surface given by the union of the two planes $014$ and $234$. Relabeling points, we may assume that $01$, $23$ and $45$ are skew lines or the 10 points are QD. 

\subsection{Ensuring the remaining four points are in general position}

We start with the assumption that the points have been labeled so that $01$, $23$ and $45$ are skew. We'd like to conclude that either the ten points are QD or there is a relabeling of all the points so that $01$, $23$, and $45$ are skew and the four remaining points $6$, $7$, $8$, and $9$ are in general position. 

There are several cases where we can easily determine if our 10 points lie on a quadric surface. 

\begin{lem} \label{lem:lines} 
(a) If the 10 points are not all distinct, then they must all lie on a quadric surface. \\
(b) If 4 of our 10 points are collinear then the 10 points lie on a quadric surface. \\
(c) More generally, if 6 of our points lie on a degree-2 plane curve, then all 10 points lie on a quadric surface.\\
(d) If 3 lines contain 9 of our 10 points then we can check whether all 10 of our points lie on a quadric surface.  \\
(e) If two lines contain six of our 10 points then we can check whether all 10 of our points lie on a quadric surface. 
\end{lem} 

\begin{proof} Parts (a) and (b) were dealt with in the introduction to this section. 
If six points lie on a degree-2 plane curve, then the six linear constraints imposed by the points are linearly dependent and all ten points lie on a quadric surface, establishing claim (c).  

Suppose three lines $L_1$, $L_2$, and $L_3$ contain 9 of our 10 points then either one of the lines contains four points and there is a quadric surface containing all 10 points by part (b), or each of the lines contains three of our points. If two of the lines meet then the ten points lie on a quadric by part (c), so we may assume that the lines are skew. Then we can use Grassmann's condition $xL_1L_2L_3x=0$ from Section \ref{ss:intro} to check whether the tenth point $x$ lies on the unique quadric surface containing the three skew lines. This establishes claim (d). 

In case (e), we may assume that the lines $L_1$ and $L_2$ are skew and each contains 3 of our points,  otherwise the 10 points lie on a quadric by parts (b) and (c). Label the remaining four points lines $a$, $b$, $c$ and $x$. By part (b) we may assume that none of these points lie on $L_1 \cup L_2$. Then there are lines through each of $a$, $b$, and $c$ hitting both $L_1$ and $L_2$ (the line $L_a$ through $a$ joins $a$ to $aL_1 \cap L_2$). 

We argue that we may assume that the three lines $L_a$, $L_b$ and $L_c$ are distinct. If $L_a = L_b$ then the linear condition on the coefficients of the quadric surface imposed by point $a$ is linearly dependent on the linear conditions imposed by the other 9 points since these nine conditions force the quadric surface to vanish on all of $L_1$ and $L_2$ and, together with the linear condition imposed by point $b$, the quadric surface must contain all three points $b$, $L_b\cap L_1$, and $L_b\cap L_2$, and hence must contain all of $L_b=L_a$. In particular, such a quadric surface must contain $a$. It follows that if $L_a=L_b$ then all ten points lie on a quadric surface; the same holds if $L_a=L_c$ or $L_b=L_c$, so the three lines $L_a$, $L_b$, and $L_c$ can be assumed distinct.  

Now we investigate the consequences if two of the three lines meet. If the line $L_a$ hits $L_b$ at point $r$, with $r \not\in L_1 \cup L_2$ then $L_1$ and $L_2$ are in the plane $abr$, violating the assumption that $L_1$ and $L_2$ are skew lines. So if $L_a$ meets $L_b$ at $r$, the point $r$ must lie on $L_1 \cup L_2$; $r \in L_1$, say. Then any quadric containing $a$, $b$, $L_1$ and $L_2$ must contain the collinear points $a$, $r$ and $L_a\cap L_2$ and so must contain all of $L_a$. Similarly, the quadric must contain all of $L_b$. But the quadric also contains the line $L_2$, so the quadric hits the plane $abr$ in three lines, $L_a$, $L_b$, and $L_2$. Any quadric meeting a plane in 3 distinct lines must contain the entire plane.
So any quadric containing $a$, $b$, $L_1$, and $L_2$ must contain the plane $abr$ and another plane containing $L_1$. So our 10 points lie in a quadric precisely when $c$ or $x$ lies in $abr$ or $c$, $x$ and $L_1$ lie in a common plane. This show that our 10 points are QD if any two of the lines $L_a$, $L_b$, and $L_c$ meet. 

So we may assume that the three lines $L_a$, $L_b$ and $L_c$ are skew. Any quadric surface containing all nine original points except (possibly) $x$, must contain $L_a$, $L_b$ and $L_c$, so it must be the unique quadric $Y$ containing all three skew lines. Then the 10 original points lie on a quadric precisely when $x$ lies on $Y$, which we can test using Grassmann's criterion, $xL_aL_bL_cx=0$. 
\end{proof} 

Recall that $01$, $23$, and $45$ are skew lines. We assume that the points $6$, $7$, $8$, and $9$ lie in a plane $\pi$ and aim to show that the ten points are QD. Set $a=01.\pi$, $b=23.\pi$, and $c=45.\pi$. Note that $a$, $b$, and $c$ are three distinct points since the lines $01$, $23$ and $45$ are skew. At most one of $a$, $b$, and $c$ can be equal to one of the points $6$, $7$, $8$ or $9$, otherwise we immediately have two skew lines that each contain 3 of our 10 points, so the 10 points are QD by Lemma \ref{lem:lines}.

The next lemma says that in certain special situations, we can permute the labels on some of the points to ensure that $01$, $23$, $45$ are skew and the remaining four points are in general position. If $\{i,j,k,\ell\}=\{6,7,8,9\}$ say that the lines $ij$ and $k\ell$ are opposite lines.  

\begin{lem}[\bf The skew swap construction] \label{lem:swap} Suppose that we choose labels $\{p,q,r,s,t,u\} = \{0,1,2,3,4,5\}$ and $\{i,j,k,\ell\} = \{6,7,8,9\}$ so that the three skew lines are $\{pq, rs, tu\}=\{01,23,45\}$, the point $pq \cap \pi$ lies off line $ij$ and both $rs\cap \pi$ and $tu \cap \pi$ lie off the opposite line $k\ell$. Swap the label of $p$ with $k$ and the label of $q$ with $\ell$. Then the three lines  $pq$, $rs$, and $tu$ are skew and the set of four remaining points $\{i,j,k,\ell\} = \{6,7,8,9\}$ lie in general position. 
\end{lem} 

\begin{proof}
Prior to making the swap, the line $pq$ cannot hit the line $ij$ since the line $pq$ hits the larger set $\pi \supset ij$ at a point off $ij$. So prior to the swap, the points $p$, $q$, $i$ and $j$ lie in general position. The same reasoning shows that prior to the swap both lines $rs$ and $tu$ miss the line $k\ell$. Since $rs$ and $tu$ were already known to be skew, the lines $k\ell$, $rs$ and $tu$ are skew. Swapping the labels on the points allows us to continue to label the skew lines $01$, $23$ and $45$ but now the remaining four points are in general position.     
\end{proof} 

In light of the skew swap construction's requirement that the points $a$, $b$, and $c$ lie off given lines, it is helpful to be able to modify our choice of skew lines to move these points, which is what the next lemma accomplishes.  

\begin{lem} \label{lem:moreskewlines} 
Suppose that the lines $ab$, $cd$, and $ef$ are skew and that $\pi$ is any plane hitting $cd$ at $g$ and hitting $ef$ at $h$. Then at least one of the following two alternatives holds: 
\begin{enumerate} 
\item the lines $ab$, $ce$, and $df$ are skew and the points $g'=ce.\pi$ and $h'=df.\pi$ are distinct and both lie off the line $gh$; 
\item the lines $ab$, $cf$, and $de$ are skew and the points $g'=cf.\pi$ and $h'=de.\pi$ are distinct and both lie off the line $gh$.  \end{enumerate} 
\end{lem}

\begin{proof} We first establish the skew properties by proving that if the lines $ab$, $ce$, and $df$ are not skew and the lines $ab$, $cf$, and $de$ are not skew, then the lines $ab$, $cd$, and $ef$ are not skew. Note that two lines $ab$ and $cd$ in 3-space meet precisely when the bracket $[abcd]$ equals zero. So three lines $ab$, $cd$, and $ef$ fail to be skew precisely when $[abcd][abef][cdef]=0$. The $4\times 4$ determinants of a $4 \times 6$ matrix are not independent, they satisfy many quadratic relations (the Grassmann-Pl\"ucker relations), which are consequences of Cramer's Rule \cite[Section 6.5]{RG}. In particular, if $[abcd] \neq 0$ then the vectors $a$, $b$, $c$, and $d$ are independent and the vector $e$ lies in their span. Cramer's Rule says that $$ \frac{[ebcd]}{[abcd]}a +
\frac{[aecd]}{[abcd]}b +
\frac{[abed]}{[abcd]}c +
\frac{[abce]}{[abcd]}d  = 
e $$ so that 
$$ [ebcd]a + [aecd]b + [abed]c + [abce]d = [abcd]e. $$ 
Wedging on the right by $a\wedge b\wedge f$ and multiplying by $[cdef]$ gives 
$$  -[abde][abcf][decf] - [abce][abdf][cedf] = [abcd][abef][cdef]. $$ 
Now both terms on the left vanish if $ab$, $de$, and $cf$ are not skew and $ab$, $ce$, and $df$ are not skew, forcing $[abcd][abef][cdef]=0$, which forces the three lines $ab$, $cd$, $ef$ to not be skew. 

In either of the two cases, $g'$ cannot equal $h'$, otherwise the two newly constructed lines meet. If $g'$ lies on $gh$ then both lines $cd$ and $ef$ lie on the plane spanned by $gh$ and $g'c$ contradicting our assumption that $cd$ and $ef$ are skew. Similarly, if $h'$ lies on $gh$ then both lines $cd$ and $ef$ lie on the plane spanned by $gh$ and $h'd$ contradicting our assumption that $cd$ and $ef$ are skew.
\end{proof}

\begin{figure}[h!t]
\begin{tikzpicture}
\tikzset{
every line/.append style ={
black,fill=black,line width=1pt}
} 
    \draw [black,fill=black] (0.7,0) circle (.05cm);
    \draw[black,fill=black] (1.2,1) circle (0.05cm);
    \draw[black,fill=black] (1.45,1.5) circle (0.05cm);
    \draw[black,fill=black] (1.7,2) circle (0.05cm);
    \draw (0.45,-0.5)--(1.95,2.5);
    \node at (1.0,0) {$6$};
    \node at (1.5,1) {$7$};
    \node at (1.75,1.5) {$8$};
    \node at (2,2) {$9$};

    \draw[black,fill=black] (3,0) circle (0.05cm);
    \draw[black,fill=black] (3.5,0) circle (0.05cm);
    \draw[black,fill=black] (4.25,0) circle (0.05cm);
    \draw (2.75,0)--(4.5,0);
    \draw[black,fill=black] (4,1.5) circle (0.05cm);
    \draw (3.4,-0.3)--(4.1,1.8);
    \draw (2.8,-0.3)--(4.2,1.8);
    \draw (4.3,-0.3)--(3.95,1.8);
    \node at (4.3,1.5) {$6$};
    \node at (2.9,0.25) {$7$};
    \node at (3.65,-0.2) {$8$};
    \node at (4.4,0.25) {$9$};

    \draw[black,fill=black] (6.1,0) circle (0.05cm);
    \draw[black,fill=black] (7.6,0) circle (0.05cm);
   \draw[black,fill=black] (6.1,0.75) circle (0.05cm);
    \draw[black,fill=black] (7.1,0.5) circle (0.05cm);
    \draw (7.9,-0.3)--(5.8,1.8);
    \draw (6.1,-0.3)--(6.1,1.8);
    \draw (5.8,0)--(9.8,0);
    \draw (5.8,0.90)--(8.1,-0.25);
    \draw (5.8,0.825)--(9.5,-0.1);
    \draw (5.8,-0.15)--(7.7,0.8);
    \node at (5.9,0.6) {$6$};
    \node at (7.2,0.75) {$7$};
    \node at (7.5,-0.2) {$8$};
    \node at (6.3,-0.2) {$9$};

\end{tikzpicture}      
\caption{Instances of cases (A) (left), (B) (middle), and (C) (right).}
\label{fig:cases}
\end{figure}
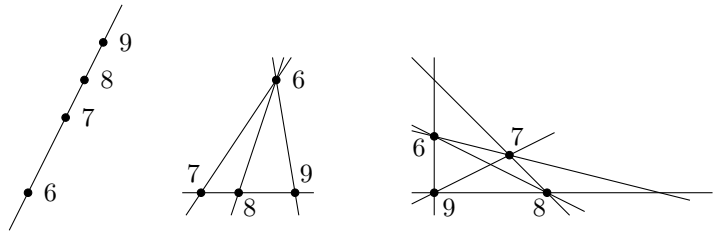

Now our argument splits into three cases (as illustrated in Figure \ref{fig:cases}): (A) all four points $6$, $7$, $8$, and $9$ are collinear; (B) three of the four points lie on a line $L$ contained in the plane $\pi$; and (C) all four points are in general position in the plane $\pi$. 

In case (A) all ten points lie on quadric by Lemma 5.1 (a). In case (B), relabel the four points so that the line $L$ contains $7$, $8$, and $9$. After relabeling, we can assume that one of the following cases holds: 
\begin{enumerate} 
\item none of the three points $a$, $b$, and $c$ lie on the union of the lines $67$, $68$, $69$, and $L$;
\item $a=6$;  
\item $a=7$; 
\item $a \in 67$ and none of $a$, $b$ or $c$ equal the four points; 
\item $a \in L$ and none of $a$, $b$ or $c$ equal any of the points $6$, $7$, $8$, or $9$. 
\end{enumerate}
In case (1), apply the skew swap construction with $a \not\in 67$ and $b \not\in 89$ and $c \not\in 89$ to produce three skew lines and four general points. In case (2), if $b$ or $c$ equals any of the four points then we are done by Lemma 5.1 (a) or (e). Then both $b$ and $c$ both lie off one of the lines $67$, $68$, and $69$; say they both lie off $69$, then apply the skew swap construction with $a \not\in 78$, $b\not\in 69$, and $c\not\in 69$ to produce three skew lines and four general points.  In case (3), if $b$ or $c$ equals any of the four points then we are done by Lemma 5.1 (e). Since $a \not\in 68$, we can apply the skew swap construction if $b$ and $c$ both lie off $L=79$. So, after relabeling we can assume $b$ lies on $L$. If $c$ is on any of the lines $67$, $68$ or $69$ then apply the skew swap construction with $c\not\in 78$, $a\not\in 69$, and $b\not\in 69$ to produce three skew lines and four general points. 
If $c$ lies on $L$ then we apply Lemma 5.2 to replace the skew lines $23$ and $45$ with new skew lines where $b$ and $c$ no longer lie on $L$. Then we have reduced to an instance of case (3) already considered.  In case (4), the points $b$ and $c$ must both lie off one of the lines $68$, $69$ or $L=89$ and we can  apply the skew swap construction (if both lie off $L=79$ then use $a \not\in 68$ in the construction). In case (5), if $b$ and $c$ both lie off $L$ then use the skew swap construction with $a \not\in 67$, $b\not\in 89$, and $c \not\in 89$. So we can assume that least one of $b$ or $c$ must lie on $L$. If $b$ lies on $L$ and $c$ lies off $L$ then use the skew swap construction with $c \not\in 89$, $a \not\in 67$, and $b\not\in 67$, while if $c$ lies on $L$ and $b$ lies off $L$ then use the skew swap construction with $b \not\in 89$, $a\not\in67$, and $c\not\in 67$. Finally, if $b$ and $c$ both lie on $L$ then use Lemma 5.2 to swap the labels on $2$ and $4$ or $2$ and $5$, to obtain new skew lines $01$, $23$, and $45$, with $b$ and $c$ no longer lying on $L$. Then we have reduced to an instance of case (5) already considered.

In case (C), after relabeling, we can assume that one of the following cases holds: 
\begin{enumerate}
\item none of the points $a$, $b$, or $c$ lie on any of the six lines $67$, $68$, $69$, $78$, $79$, and $89$; 
\item $a = 6$; 
\item $a$ lies on $67$ and none of $a$, $b$, or $c$ equals any of the four points;
\item $a = 67 \cap 89$ and neither $b$ nor $c$ equals any of the points $6$, $7$, $8$, or $9$. 
\end{enumerate}
Cases (1) and (2) are handled exactly as in cases (1) and (2) in case (B) above. In case (3), both $b$ and $c$ lie off one of the five lines $67,68,69,78,79$ and then this line, $23$, 
and $45$ are skew and the remaining four of our 10 points are in general position. In case (4), none of $a$, $b$, or $c$ equals $6$, $7$, $8$, or $9$, otherwise we are in case (2). If both $b$ and $c$ lie off one of the lines $68$, $69$, $78$, or $79$ and then this line, $23$, 
and $45$ are skew and the remaining four of our 10 points are in general position. So one of $b$ or $c$ must lie on each of the lines $68$, $69$, $78$, or $79$. When combined with $\{b,c\} \cap \{6,7,8,9\} = \emptyset$, this forces 
$\{b,c\} = \{ 68\cap 79, 69\cap 78\}$. Then we can use Lemma 5.2 to replace the skew lines $23$ and $45$ with new skew lines where $b$ and $c$ no longer lie at $68\cap 79$ and $69\cap 78$. Then we have reduced to an instance already considered.

The three cases, (A), (B), and (C) exhaust all possibilities, so we conclude that either the 10 points are QD or, after relabeling, the lines $01$, $23$, and $45$ are skew and the points $6$, $7$, $8$, and $9$ are in general position. This is precisely the case treated in Section \ref{ss:generic}, where we show that the 10 points are QD. So, it turns out that the 10 points are always quadric decidable. We have a series of synthetic constructions that can be used to determine if the 10 points lie on a quadric surface. 

\section{Open Problems and Further Reading}

The synthetic methods developed here could be applied to many related problems. One obvious place to start is to use the methods to recover Pascal's Theorem. It would be good to understand how to reduce the resulting construction to the simple statement in Theorem \ref{thm:PBM}. Recently, Traves and Wehlau \cite{TW} devised a synthetic method to check whether 10 points lie on a plane cubic curve. The construction finds two degree-4 curves, each the union of two conics, whose intersection contains the 10 points. The Cayley-Bacharach Theorem \cite{EGH} implies that 10 points lie on a cubic precisely when the residual set of 6 points lies on a conic. Unfortunately, it is not possible to precisely locate all six residual points using a purely synthetic construction; however, synthetic methods still suffice to check whether the six points lie on a conic (see the article for details). Traves \cite{T} also describes how Grassmann's methods can be used to determine whether 10 points lie on a plane cubic, as well as other synthetic geometry problems involving cubics and conics. 

Can our methods be used to produce a synthetic construction that checks when 15 points in the plane lie on a degree-4 curve? In this direction, we note that if the 15 points lie among the 25 points of intersection of two degree-5 curves, each the product of five linear forms, then, using the Cayley-Bacharach Theorem,  the problem immediately reduces to checking that the 10 residual points lie on a cubic curve, for which we can use the synthetic methods described above. 

Caminata and Schaffler \cite{CS} gave synthetic conditions for $d+4$ points to lie on a rational normal curve of degree $d$ in $\mathbb{P}^4$. Ciliberto and Miranda \cite[Theorem 3]{CM} simplified the conditions: the points lie on a rational normal curve precisely when both 
\begin{enumerate}
\item the projection from one of the points $p_i$ to a $\mathbb{P}^{d-1}$ applied to the remaining points give $d+3$ points on a rational normal curve, and 
\item the projection from the span of $d-2$ points (distinct from $p_i$) to a $\mathbb{P}^2$ applied to the remaining points gives six points lying on a conic. 
\end{enumerate} 
Ciliberto and Miranda discuss several generalizations of Pascal's Theorem, investigating particular instances under which ten points lie on a quadric surface in $\mathbb{P}^3$. For example, they note that if the ten points are broken into two groups of five, then there are rational normal curves through each group and these meet in 5 points precisely when the ten points lie on a quadric surface \cite[Proposition 7]{CM}. 

Ciliberto and Miranda describe another remarkable result, originally due to Richmond, but rediscovered\footnote{It is distressing that so much classical algebraic geometry gets forgotten, but it is good to see these results restated and proved using modern methods.} by both Segre and Brown: five ordered lines $L_1,\ldots,L_5$ lie on a quadric surface in $\mathbb{P}^4$ precisely when the five points $R_i = L_i\wedge(L_{i-1}L_{i+1})$ (indicies taken modulo 5) lie on a hyperplane.  

Generalizing by degree rather than dimension, we ask for a synthetic construction that checks when 20 points in 3-space lie on a cubic surface. It would be best if the solutions to synthetic problems posed in $\mathbb{P}^n$ only require synthetic constructions in $\mathbb{P}^n$ or projective spaces of lower dimension, though perhaps this is too much to require. 

An underlying idea in the paper is the identification of $\mathbb{P}^3$ with a tetrahedron. This suggests that we can ask for synthetic proofs of geometric results on other toric varieties. What would such constructions look like? 

Synthetic geometry results can all be stated as incidence theorems. Fomin and Pylyavskyy \cite{FP} give a construction involving quadrilateral tilings of a closed oriented surface 
that produces new incidence theorems and generalizes many known results, including the theorems of Pappus, Desargues, and M\"{o}bius. Richter-Gebert \cite{RG-Coxeter} investigated similar constructions using the theorems of Ceva and Menelaus as  building blocks for incidence results.  It would be fascinating to know if all incidence theorems arise from their construction(s), or even to find examples that cannot be produced using these constructions. 
 
\smallskip
\noindent {\bf Acknowledgements}:  J\"{u}rgen Richter-Gebert and Susanne Apel were incredibly gracious and met with me in Munich to discuss mathematics that eventually led to this paper. This work also benefited tremendously from conversations with Mike Roth, David Wehlau, Chris Peterson, and Rick Miranda. I'm grateful for their friendship and advice. 

The main results of this article were originally presented on a poster at the Combinatorial, Computational, and Applied Algebraic Geometry conference in Seattle during the Summer of 2022. The conference celebrated Bernd Sturmfels and his mathematical contributions on the occasion of his 60th birthday. I thank the organizers for the opportunity to present my work there.   

\bibliographystyle{plain}
\bibliography{sample}

\begin{thebibliography}{10}

\bibitem{SA}
Susanne Apel.
\newblock Cayley factorization and the area principle.
\newblock {\em Discrete Comput. Geom.}, 55(1):203--227, 2016.

\bibitem{JRGSA}
Susanne Apel and J\"urgen Richter-Gebert.
\newblock Cancellation patterns in automatic geometric theorem proving.
\newblock In {\em Automated deduction in geometry}, volume 6877 of {\em Lecture
  Notes in Comput. Sci.}, pages 1--33. Springer, Heidelberg, 2011.

\bibitem{CS}
Alessio Caminata and Luca Schaffler.
\newblock A {P}ascal's theorem for rational normal curves.
\newblock {\em Bull. Lond. Math. Soc.}, 53(5):1470--1485, 2021.

\bibitem{CM}
Ciro Ciliberto and Rick Miranda.
\newblock Variations on pascal's theorem, 2024.

\bibitem{CG}
H.~S.~M. Coxeter and S.~L. Greitzer.
\newblock {\em Geometry revisited}, volume~19 of {\em New Mathematical
  Library}.
\newblock Random House, Inc., New York, 1967.

\bibitem{EGH}
David Eisenbud, Mark Green, and Joe Harris.
\newblock Cayley-{B}acharach theorems and conjectures.
\newblock {\em Bull. Amer. Math. Soc. (N.S.)}, 33(3):295--324, 1996.

\bibitem{FP}
Sergey Fomin and Pavlo Pylyavskyy.
\newblock Incidences and tilings, 2023.

\bibitem{HG}
Hermann Grassmann.
\newblock {\em Extension theory}, volume~19 of {\em History of Mathematics}.
\newblock American Mathematical Society, Providence, RI; London Mathematical
  Society, London, 2000.
\newblock Translated from the 1896 German original and with a foreword,
  editorial notes and supplementary notes by Lloyd C. Kannenberg.

\bibitem{M2}
Daniel~R. Grayson and Michael~E. Stillman.
\newblock Macaulay2, a software system for research in algebraic geometry.

\bibitem{RG-Coxeter}
J\"urgen Richter-Gebert.
\newblock Meditations on {C}eva's theorem.
\newblock In {\em The {C}oxeter legacy}, pages 227--254. Amer. Math. Soc.,
  Providence, RI, 2006.

\bibitem{RG}
J\"{u}rgen Richter-Gebert.
\newblock {\em Perspectives on projective geometry}.
\newblock Springer, Heidelberg, 2011.
\newblock A guided tour through real and complex geometry.

\bibitem{S-AIT}
Bernd Sturmfels.
\newblock {\em Algorithms in invariant theory}.
\newblock Texts and Monographs in Symbolic Computation. SpringerWienNewYork,
  Vienna, second edition, 2008.

\bibitem{SW}
Bernd Sturmfels and Walter Whiteley.
\newblock On the synthetic factorization of projectively invariant polynomials.
\newblock {\em J. Symbolic Comput.}, 11(5-6):439--453, 1991.
\newblock Invariant-theoretic algorithms in geometry (Minneapolis, MN, 1987).

\bibitem{T}
Will Traves.
\newblock What grassmann knew: Incidence theorems on cubics, 2023.

\bibitem{TW}
Will Traves and David Wehlau.
\newblock Ten points on a cubic.
\newblock {\em Amer. Math. Monthly}, 131(2):112--130, 2024.

\bibitem{TY}
H.W. Turnbull and Alfred Young.
\newblock The linear invariants of ten quaternary quadrics.
\newblock {\em Transactions of the Cambridge Philosophical Society},
  23(10):264--302, 1927.

\bibitem{White90}
Neil White.
\newblock Implementation of the straightening algorithm of classical invariant
  theory.
\newblock In {\em Invariant theory and tableaux ({M}inneapolis, {MN}, 1988)},
  volume~19 of {\em IMA Vol. Math. Appl.}, pages 36--45. Springer, New York,
  1990.

\end{thebibliography}

\end{document}